\documentclass{article}

\usepackage{
 amsmath,
 amsxtra,
 amsthm,
 amssymb,
 etex,
 mathrsfs,
 stmaryrd
  }
\usepackage[all]{xy}
\usepackage{hyperref}

\newtheorem{theorem}{Theorem}[section]
\newtheorem{lemma}[theorem]{Lemma}

\newtheorem{proposition}[theorem]{Proposition}
\newtheorem{corollary}[theorem]{Corollary}
\newtheorem*{thm2}{Theorem}

\theoremstyle{definition}
\newtheorem{defn}[theorem]{Definition}
\newtheorem{remark}[theorem]{Remark}

\newcommand{\bd}{\begin{defn}}
\newcommand{\ed}{\end{defn}}
\newcommand{\bl}{\begin{lemma}}
\newcommand{\el}{\end{lemma}}
\newcommand{\bp}{\begin{proposition}}
\newcommand{\ep}{\end{proposition}}
\newcommand{\bt}{\begin{theorem}}
\newcommand{\et}{\end{theorem}}
\newcommand{\bc}{\begin{corollary}}
\newcommand{\ec}{\end{corollary}}
\newcommand{\br}{\begin{remark}}
\newcommand{\er}{\end{remark}}
\newcommand{\ba}{\begin{array}}
\newcommand{\ea}{\end{array}}
\newcommand{\bpf}{\begin{proof}}
\newcommand{\epf}{\end{proof}}

\newcommand{\Q}{\mathbb{Q}}

\newcommand{\Z}{\mathbb{Z}}
\newcommand{\Zp}{\mathbb{Z}_p}
\newcommand{\Op}{\mathcal{O}}
\newcommand{\Ga}{\Gamma}
\newcommand{\ga}{\gamma}

\newcommand{\e}{\varepsilon}
\newcommand{\ze}{\zeta}
\newcommand{\om}{\omega}

\DeclareMathOperator{\Gal}{Gal}

\DeclareMathOperator{\rank}{rank}

\DeclareMathOperator{\Cl}{Cl}

\newcommand{\Iw}{\mathrm{Iw}}

\newcommand{\Tr}{\mathrm{Tr}}

\newcommand{\lra}{\longrightarrow}

\newcommand{\ot}{\otimes}

\newcommand{\cyc}{\mathrm{cyc}}
\newcommand{\ps}[1]{\llbracket #1 \rrbracket}

\newcommand{\ilim}{\displaystyle \mathop{\varinjlim}\limits}
\newcommand{\plim}{\displaystyle \mathop{\varprojlim}\limits}

\numberwithin{equation}{section}

\begin{document}
\title{On capitulations of even $K$-groups and pseudo-null submodules
in $\Zp^d$-extensions}
 \author{
  Meng Fai Lim\footnote{School of Mathematics and Statistics, Key Laboratory of Nonlinear Analysis and Applications (Ministry of Education),
Central China Normal University, Wuhan, 430079, P.\ R.\ China.
 E-mail: \texttt{limmf@ccnu.edu.cn}} }
\date{}
\maketitle

\begin{abstract} \footnotesize
\noindent  The capitulations of
ideals in $\Zp^d$-extensions and pseudo-null submodules of the classical Iwasawa modules
are closely related as evidenced in the works of Ozaki and Fujii. In this paper, we investigate the analogous situation for the even $K$-groups. As an application, we obtain a new sufficient condition for existence of non-trivial pseudo-null submodules in the classical Iwasawa modules.


\end{abstract}

\section{Introduction} \label{intro}

In this paper, $p$ will always denote an odd prime. Let $F$ be a number field and let $K$ be a $\Zp$-extension of $F$. For every finite extension $L$ of $F$ contained in $K$, we let $A_L$ denote the $p$-primary part of the ideal class group of $\Op_L$, where $\Op_L$ is the ring of integer of $L$. For an extension $L'$ of $L$ contained in $K$, there is a natural map $A_L \lra A_{L'}$ induced by the inclusion $\Op_L\subseteq \Op_{L'}$. On the other hand, we have a map $A_{L'} \lra A_{L}$
going the other way which is induced by the (ideal) norm. Then one has
\[  \ilim_L A_L \quad \mbox{and}\quad \plim_L A_L,\]
where the direct limit (resp., the inverse limit) is taken with respect to the inclusion maps (resp., the norm maps). These two limit modules come naturally equipped with  $\Zp\ps{\Ga}$-module structures, where $\Ga= \Gal(K/F)\cong\Zp$. Write $A_K = \ilim_L A_L$
and $X_K = \plim_L A_L$. Let $X_K^0$ be the maximal finite submodule of $X_K$.
The following result of Ozaki \cite{Oz} serves as the main motivation of our paper.

\begin{thm2}[Ozaki]
Suppose that every ramified prime of $K/F$ is totally ramified. Then $X^{0}_K\neq 0$ if and only if $\ker\left(A_L \lra A_K\right)\neq 0$ for some finite extension $L$ of $F$ contained in $K$.
\end{thm2}

When $F$ is a totally real field and $K$ is the cyclotomic $\Zp$-extension of $F$, the non-triviality of $X_K^0$ corresponds to a weakened form of Greenberg's conjecture (for instances, see \cite{Fuj20, NQD, NQD2}; for the original Greenberg's conjecture, we refer readers to the work of Greenberg \cite{Gr76}).

Let $\widetilde{F}$ denote the compositum of all $\Zp$-extensions of $F$. We then define the module $X_{\widetilde{F}}$ in a similar manner by taking the inverse limit of $A_L$ with $L$ running over all finite subextensions of $\widetilde{F}/F$. In this context, one central theme of study is the nontriviality of $X_{\widetilde{F}}^0$, the maximal pseudo-null $\Zp\ps{\Gal(\widetilde{F}/F)}$-submodule of $X_{\widetilde{F}}$. This latter property may be regarded as a weakened form of Greenberg's generalized conjecture. (For a selection of related investigations exploring this direction of research, we refer the reader to the works \cite{Fuj10, Kl, Mur}, where we emphasize that this list is far from being exhaustive. For the statement of Greenberg's generalized conjecture, readers may consult \cite[Conjecture 3.5]{Gr98}.) Turning to this topic, a recent generalization of Ozaki's result to the setting of $\Zp^d$-extensions with $d\geq 2$ has been established by Fujii \cite{Fuj}, whose results we will present below.

\begin{thm2}[Fujii]
Let $F_\infty$ be a $\Zp^d$-extension of $F$.
Suppose that $F$ has only one prime above $p$ and this prime is totally ramified in $F_\infty/F$. If $X^{0}_{F_\infty}\neq 0$, then $\ker\left(A_L \lra A_{F_\infty}\right)\neq 0$ for some finite extension $L$ of $F$ contained in $F_\infty$.
\end{thm2}

\begin{thm2}[Fujii]
Let $F_\infty$ be a $\Zp^d$-extension of $F$.
Suppose that $F$ has only one prime above $p$ and this prime is totally ramified in $F_\infty/F$. Assume further that $A_F$ is a non-trivial cyclic group. Then $\ker\left(A_L \lra A_K\right)\neq 0$ for some finite extension $L$ of $F$ contained in $F_\infty$, then $X^{0}_{F_\infty}\neq 0$.
\end{thm2}

In this paper, we will be interested in the natural analogous situation for higher even $K$-groups. Consistent with our earlier conventions, we continue to let $p$ denote an odd prime, and let $F$ be a number field. For a ring $R$ with identity, write $K_n(R)$ for the algebraic $K$-groups of $R$ in the sense of Quillen \cite{Qui73a}, and we refer readers to \cite{Kol, WeiKbook} for additional standard references on these groups. Let $\Op_F$ denote the ring of integers of $F$. The foundational results of Garland \cite{Gar}, Quillen \cite{Qui73b} and Borel \cite{Bo} collectively establish that the group $K_{2i}(\Op_F)$ is finite for every $i\geq 1$.

Now let $F_\infty$ be a $\Zp^d$-extension of $F$. For any pair of finite subextensions with $L\subseteq L'$, the natural inclusion $\Op_L\lra \Op_{L'}$ of rings induces a map
\[\jmath_{L/L'}: K_{2i}(\Op_L)[p^\infty]\lra K_{2i}(\Op_{L'})[p^\infty]\]
 by the functoriality of algebraic $K$-theory. We should mention that this induced map may not be injective in general. In the reverse direction, the norm map (also called the transfer map) gives a homomorphism
 \[\Tr_{L'/L}: K_{2i}(\Op_{L'})[p^\infty]\lra K_{2i}(\Op_{L})[p^\infty].\]
  With these transition maps in hand, we define the direct limit and inverse limit
\[ \ilim_L K_{2i}(\Op_L)[p^\infty]\quad\mbox{and}\quad \plim_L K_{2i}(\Op_L)[p^\infty],\]
where the direct limit uses $\jmath_{L/L'}$ and inverse limit uses the norm maps $\Tr_{L'/L}$. Both limit modules come equipped with natural $\Zp\ps{G}$-module structures, where $G=\Gal(F_\infty/F)\cong \Zp^d$. We shall denote the direct limit module by $A_{2i,F_\infty}$ and the inverse limit module by $X_{2i,F_\infty}$. Finally, we denote by $X^{0}_{2i,F_\infty}$ the maximal pseudo-null submodule of $X_{2i,F_\infty}$.

Our first main result, which serves as the direct analogue of Ozaki's result in the setting of even $K$-groups, is as follows.

\bt[Theorem \ref{mainthm1}]
Suppose that $d=1$, i.e., $F_\infty$ is a $\Zp$-extension of $F$. Then $X^{0}_{2i,F_\infty}\neq 0$ if and only if there exists a finite extension $L$ of $F$ contained in $F_\infty$ such that
\[ \ker\Big(K_{2i}(\Op_L)[p^\infty] \lra A_{2i,F_\infty}\Big)\neq 0.\]
\et

The subsequent two results focus on the setting of multiple $\Zp$-extensions, and are the exact even $K$-group counterparts of \cite[Theorems 1.2 and 1.3]{Fuj} from Fujii's prior work.

\bt[Theorem \ref{mainthm2a}]
Let $F_\infty$ be a $\Zp^d$-extension of $F$ with $d\geq 2$. Suppose that $X^{0}_{2i,F_\infty}\neq 0$. Then there exists a finite extension $L$ of $F$ contained in $F_\infty$ such that
\[ \ker\Big(K_{2i}(\Op_L)[p^\infty] \lra A_{2i, F_\infty}\Big)\neq 0.\]
\et

\bt[Theorem \ref{mainthm2b}]
Let $F_\infty$ be a $\Zp^d$-extension of $F$ with $d\geq 2$. Suppose that $K_{2i}(\Op_F)[p^\infty]$ is a nontrivial cyclic group and there exists a finite extension $L$ of $F$ contained in $F_\infty$ such that
\[ \ker\Big(K_{2i}(\Op_L)[p^\infty] \lra A_{2i, F_\infty}\Big)\neq 0.\]
Then $X^{0}_{2i,F_\infty}\neq 0$.
\et

As a byproduct of our proof for the above theorem, we obtain the following result, which not only acts as a key intermediate step for the main argument, but may also be of independent interest (see Remark \ref{H Iw remark}(2)).

\bt[Theorem \ref{projdim}]
Let $F_\infty$ be a $\Zp^d$-extension of $F$ with $d\geq 2$. Suppose that there exists a finite extension $L$ of $F$ contained in $F_\infty$ with
\[ \ker\Big(K_{2i}(\Op_L)[p^\infty] \lra A_{2i, F_\infty}\Big)\neq 0.\]
Then the projective dimension of $X_{2i,F_\infty}$ as a $\Zp\ps{G}$-module is at least $2$.
\et

Our proofs are modelled after the works of Ozaki \cite{Oz, Oz01} and Fujii \cite{Fuj}. Unlike class groups, the even $K$-groups enjoy stronger Galois descent properties, as recorded in Proposition \ref{gal descent}. This allows us to eliminate the strict constraints on the number of $p$-adic primes of $F$ and their ramification that are required in Fujii's theorems. It is therefore a natural line of inquiry to ask whether our $K$-theoretic results can yield new complementary statements to those of Fujii's. This forms the natural next subject that we will now describe.

Fix once and for all a primitive $p$-th root of unity $\zeta\in\overline{\Q}$, a generator $g$ of $(\Z/p\Z)^{\times}$, and let $\sigma\in\Gal(\Q(\mu_p)/\Q)$ denote the map $\zeta\mapsto\zeta^g$. Let $r=|\Gal(F(\mu_p)/F)|$, and so the elements of the group $\Gal(F(\mu_p)/F)$ may be identified with elements of the form $\sigma^{k(p-1)/r}$, where $k=0,1,..., r-1$. For a given integer $j$, we define
\[ \e_j:=\e_j(F): = \frac{1}{r}\sum_{k=0}^{r-1}\omega(g)^{jk(p-1)/r}\sigma^{-k(p-1)/r}\in \Zp[\Gal(F(\mu_{p})/F)],\]
where $\omega$ is the Teichm\"{u}ller character of the group $(\Z/p\Z)^{\times}$. We now introduce the following property.
\begin{itemize}
  \item[$(\mathbf{Sp})$] Every prime of $F(\ze+\ze^{-1})$ above $p$ does not split in $F(\mu_p)$.
\end{itemize}
Finally, we denote by $F^\cyc$ the cyclotomic $\Zp$-extension of $F$. With these notations in place, we can state the following.

\bt[Theorem \ref{class group}]
Let $F$ be a number field which satisfies $(\mathbf{Sp})$, and let $F_\infty$ be a $\Zp^d$-extension of $F$ such that $F^\cyc\subseteq F_\infty$. Suppose that there exists an odd integer $i$ with $1\leq i\leq r-1$ that satisfies all of the following conditions.
 \begin{enumerate}
   \item[$(a)$] There exists a finite extension $L$ of $F$ contained in $F_\infty$ such that
\[ \ker\Big(K_{2i}(\Op_{L})[p^\infty] \lra A_{2i,F_\infty}\Big)\neq 0.\]
\item[$(b)$] $K_{2i}(\Op_F)[p^\infty]$ is a non-trivial cyclic group.
  \end{enumerate}
Then $X_{E_\infty}^0\neq 0$, where $E_\infty = F_\infty(\mu_p)$.
\et

We end the introductory section with a brief overview of the paper. Section \ref{K and Galois} collects required background results on even $K$-groups and their relation to Galois cohomology. Our main $K$-theoretic results are proved in Section \ref{sec:main}, and Section \ref{sec:class} applies these results to give a sufficient condition for the existence of nontrivial pseudo-null submodules.

\subsection*{Acknowledgments}
The author is partially supported
by the Open Research Fund of Hubei Key Laboratory of Mathematical Sciences (Central China Normal University) No. MPL2025ORG020.

\section{$K$-groups and Galois cohomology} \label{K and Galois}

In the present section, we recall key known facts about even $K$-groups and their relation to Galois cohomology, while retaining all notations defined in Section \ref{intro}. Additional necessary notations will also be introduced as we proceed.

As a start, let $\mu_{p^n}$ denote the group of all $p^n$-th power root of unity, and let $S$ be the set of all primes of $F$ consisting precisely of the infinite primes and the primes above $p$. Denote by $F_S$ the maximal algebraic extension of $F$ unramified outside $S$. For every extension $\mathcal{L}$ of $F$ contained in $F_S$, we write $G_S(\mathcal{L})$ for the Galois group $\Gal(F_S/\mathcal{L})$.

The natural action of $G_S(F)$ on $\mu_{p^\infty}$ induces the cyclotomic character
\[\chi: G_S(F) \lra \mathrm{Aut}(\mu_{p^\infty}) \cong \Zp^{\times}.\]
For a $\Zp[G_S(F)]$-module $X$, we denote by $X(i)$ the $i$-fold Tate twist of $X$. More precisely, $X(i)$ is the $G_S(F)$-module which is $X$ as a $\Zp$-module but with a $G_S(F)$-action redefined by the rule
\[ \sigma\cdot x = \chi(\sigma)^i\sigma x,  \]
where the action on the right is the original action of $G_S(F)$ on $X$.

In his foundational work \cite{Sou}, Soul\'e established a connection between the higher $K$-groups with (continuous) Galois cohomology groups via the $p$-adic Chern class maps
\[ \mathrm{ch}_{i,k}^{(p)}: K_{2i+2-k}(\Op_{F,S})\ot \Zp \lra H^k(G_{S}(F), \Zp(i+1))\]
for $i\geq 1$ and $k =1,2$. (We refer the reader to loc.\ cit.\ for the full, precise construction of these maps.) Thanks to the deep work of Rost and Voevodsky \cite{Vo} (also see \cite{Wei09}), these $p$-adic Chern class maps are now known to be isomorphisms. As a direct consequence of this isomorphism result, we obtain the following identification between the Sylow $p$-subgroups of even $K$-groups and second Galois cohomology groups.

\bp \label{K2 = H2}
For $i\geq 1$, one has the following isomorphisms
\[ K_{2i}(\Op_{F})[p^\infty] \cong H^2\big(G_{S}(F), \Zp(i+1)\big).\]
\ep

Furthermore, from the discussion in \cite[Chap.\ III]{Sou}, the aforementioned isomorphism fits into the following commutative diagrams
\begin{equation} \label{N cor}\entrymodifiers={!! <0pt, .8ex>+} \SelectTips{eu}{}
\xymatrix{
      K_{2i}(\Op_{L,S})[p^\infty] \ar[r]^(.45){\mathrm{ch}_i^L} \ar[d]^{\Tr_{L/F}} &  H^2\left(G_S(L), \Zp(i+1)\right) \ar[d]^{\mathrm{cor}} \\
      K_{2i}(\Op_{F,S})[p^\infty] \ar[r]^(.45){\mathrm{ch}_i^F} &  H^2\left(G_S(F), \Zp(i+1)\right)}
      \end{equation}
      \begin{equation} \label{j tr}\entrymodifiers={!! <0pt, .8ex>+} \SelectTips{eu}{}
\xymatrix{
      K_{2i}(\Op_{F,S})[p^\infty] \ar[r]^(.45){\mathrm{ch}_i^L} \ar[d]^{j_{L/F}} &  H^2\left(G_S(F), \Zp(i+1)\right) \ar[d]^{\mathrm{res}} \\
      K_{2i}(\Op_{L,S})[p^\infty] \ar[r]^(.45){\mathrm{ch}_i^F} &  H^2\left(G_S(L), \Zp(i+1)\right)}
      \end{equation}

We now introduce the Iwasawa cohomology groups which will be an important object for our study. For every extension $\mathcal{L}$ of $F$ contained in $F_S$, the Iwasawa cohomology groups are defined by
 \[H^k_{\Iw,S}\big(\mathcal{L}/F, \Zp(i+1)\big):= \plim_L H^k\big(G_S(L), \Zp(i+1)\big), \]
 where $L$ runs through all finite extensions of $F$ contained in $\mathcal{L}$ and the transition maps are given by the corestriction maps. To simplify notation,
we shall omit the subscript `$S$' from the Iwasawa cohomology groups in all subsequent discussions. By \cite[Theorem 1]{Iw73}, if $F_\infty$ is a $\Zp^d$-extension of $F$, we then have $F_\infty\subseteq F_S$. Taking the commutative diagram (\ref{N cor}) into account, we have the identification
\[  X_{2i,F_\infty}:=\plim_{L} K_{2i}(\Op_{L})[p^\infty]\cong H^2_{\Iw}\big(F_\infty/F, \Zp(i+1)\big)\]
 that will be used in all subsequent discussions without further mention.

The groups $H^k_{\Iw}\big(F_\infty/F, \Zp(i+1)\big)$ come naturally equipped with $\Zp\ps{G}$-module structures and can be shown to be finitely generated over $\Zp\ps{G}$ (see \cite[Proposition 4.1.3]{LimSh}). For the second Iwasawa cohomology group, which is the main object of interest, an even stronger property holds, as we state in the following result.

\bp \label{tor}
Let $i\geq 1$. For every $\Zp^d$-extension $F_\infty$ of $F$, the module $H^2_{\Iw}\big(F_\infty/F, \Zp(i+1)\big)$ is torsion over $\Zp\ps{\Gal(F_\infty/F)}$.
\ep

\bpf
This is a special case of \cite[Proposition 4.1.1]{LimKgroups}, where the torsionness assertion is in fact established for general, possibly non-commutative, $p$-adic Lie extensions.
\epf

The following version of Tate's descent spectral sequence for Iwasawa cohomology will be crucial for our subsequent discussion.

\bp \label{gal descent}
Let $F_\infty$ be $\Zp^d$-extension of $F$ and assume $i\geq 1$.
Let $U$ be a closed normal subgroup of $G=\Gal(F_\infty/F)$, and write $\mathcal{L}$ for the fixed field of $U$ in $F_\infty$. Then we have a homological spectral sequence
 $$ H_r\big(U, H^{2-s}_{\Iw}(F_\infty/F, \Zp(i+1))\big)\Longrightarrow H^{2-r-s}_{\Iw}\big(\mathcal{L}/F, \Zp(i+1)\big) $$
and an isomorphism
\[ H^2_{\Iw}\big(F_\infty/F, \Zp(i+1)\big)_U \cong H^2_{\Iw}\big(\mathcal{L}/F, \Zp(i+1)\big).\]
Furthermore, if $U$ is taken to be an open subgroup of $G$, then its fixed field $\mathcal{L}$ is a finite extension of $F$ and we have the isomorphism
\[ H^2_{\Iw}\big(F_\infty/F, \Zp(i+1)\big)_U \cong H^2\big(G_S(\mathcal{L}), \Zp(i+1)\big)\cong K_{2i}(\Op_{\mathcal{L}})[p^\infty].\]
\ep

\bpf
The spectral sequence follows from a result of Nekov\'a\v{r} \cite[Proposition 4.2.3]{Ne}. (One can also find more general versions that imply these in \cite[Proposition 1.6.5]{FK} or \cite[Theorem 3.1.8]{LimSh}.) The isomorphisms of the proposition follow from reading off the initial term of the spectral sequence. Finally, in the event that $U$ is an open subgroup of $G$, then $\mathcal{L}$ is a finite extension of $F$, and so $H^2_\Iw(L/F,\Zp(i+1))$ identifies with $H^2(G_S(L), \Zp(i+1))$.
\epf

\section{Main results} \label{sec:main}

We are now in a position to prove our main results.  We begin with the following statement, which is the direct analogue of Ozaki's result in the context of even $K$-groups.

\bt \label{mainthm1}
Let $F_\infty$ be a $\Zp$-extension of $F$. Then we have $X^{0}_{2i,F_\infty}\neq 0$ if and only if
\[ \ker\Big(K_{2i}(\Op_L)[p^\infty] \lra A_{2i,F_\infty}\Big)\neq 0.\]
\et

\bpf
For the proof, we shall write $\Ga$ for the Galois group $\Gal(F_\infty/F)$. We then identify $\Zp\ps{\Ga}$ with the power series ring $\Zp\ps{T}$ in one variable. Set $w_n = (T+1)^{p^n}-1$ and $v_{m,n}= w_m/w_n$ for $m>n\geq 0$. Let $F_n$ denote the unique intermediate extension of $F_\infty/F$ with $[F_n : F] = p^n$, and set $\Gamma_n = \Gal(F_\infty/F_n)$. We claim that
\begin{equation}\label{kercap}
 \ker\Big(K_{2i}(\Op_{F_n})\lra A_{2i,F_\infty}\Big) = \mathrm{im}\Big(X_{2i, F_\infty}^0\lra X_{2i, F_\infty}/w_n = K_{2i}(\Op_{F_n})\Big).
\end{equation}
Assuming the claim for the moment, the implication $(\Leftarrow)$ of the main theorem follows directly. Conversely, if $X_{2i}^0\neq 0$, then since $\bigcap_n w_nX_{2i,F_\infty}=0$, we must have $X_{2i}^0\nsubseteq w_nX_{2i,F_\infty}$ for some $n$. By our claim, this in turn implies that $\ker\Big(K_{2i}(\Op_L)[p^\infty] \lra A_{2i,F_\infty}\Big)\neq 0$ as required.

It therefore remains to verify our claim (\ref{kercap}). Since $X_{2i,F_\infty}^0$ is finite, we can find $m>n$ such that $w_m/w_n$ annihilates $X_{2i,F_\infty}^0$. Consider the following diagram
\[ \entrymodifiers={!! <0pt, .8ex>+} \SelectTips{eu}{}\xymatrix{
    X_{2i,F_\infty}/w_n \ar[r]_(.4){\sim}\ar[d]_{w_m/w_n} &    H^2\big(G_S(F_n), \Zp(i+1)\big) \ar[d]_{\mathrm{res}} & \ar[l]^(.4){\sim} K_{2i}(\Op_{F_n})[p^\infty] \ar[d]_{j_{n,m}}\\
    X_{2i,F_\infty}/w_m \ar[r]_(.4){\sim} &  H^2\big(G_S(F_m), \Zp(i+1)\big) &  \ar[l]^(.4){\sim}  K_{2i}(\Op_{F_m})[p^\infty]
     }\]
The commutativity of the left square follows from \cite[Lemma 3.6]{LimKlimit}, while the right square commutes by (\ref{j tr}). The ``$\supseteq$'' part of claim (\ref{kercap}) is an immediate consequence of the diagram. For the reverse inclusion, we let $x\in \ker\Big(K_{2i}(\Op_{F_n})[p^\infty]\lra A_{2i,F_\infty}\Big)$. Then $x\in \ker\Big(K_{2i}(\Op_{F_n})[p^\infty]\lra K_{2i}(\Op_{F_m})[p^\infty]\Big)$ for some big enough $m>n$. Choose a lift $\widetilde{x} \in X_{2i,F_\infty}$ of $x$ via the surjection $X_{2i,F_\infty} \twoheadrightarrow X_{2i,F_\infty}/w_n \cong K_{2i}(\mathcal{O}_{F_n})[p^\infty]$. By the commutative diagram, there exists  $y\in X_{2i,F_\infty}$ such that $\frac{w_m}{w_n}\widetilde{x} = w_my$, or equivalently, $\frac{w_m}{w_n}(\widetilde{x}-w_n y)=0$. As $X_{2i,F_\infty}/w_m\cong K_{2i}(\Op_{F_m})[p^\infty]$ is finite, the group $X_{2i,F_\infty}[w_m]$ is finite, and hence $X_{2i,F_\infty}[v_{m,n}]$ is also finite. In particular, $X_{2i,F_\infty}[w_m/w_n]\subseteq X_{2i,F_\infty}^0$. Thus, $\widetilde{x}-w_n y\in X_{2i,F_\infty}^0$. Since the image of this element under the map $X_{2i, F_\infty}^0\lra X_{2i, F_\infty}/w_n = K_{2i}(\Op_{F_n})$ is $x$, we have the reverse inclusion.
\epf

We now come to the algebraic $K$-theoretical analogue of \cite[Theorem 1.3]{Fuj}.

\bt \label{mainthm2a}
Suppose that $F_\infty$ is a $\Zp^d$-extension of $F$ with $d\geq 2$. If $X^{0}_{2i,F_\infty}\neq 0$, then there exists a finite extension $L$ of $F$ contained in $F_\infty$ such that
\[ \ker\Big(K_{2i}(\Op_L)[p^\infty] \lra A_{2i, F_\infty}\Big)\neq 0.\]
\et

Before giving the proof of this said theorem, we first establish a few preparatory lemmas. The   following lemma is the first of these.

\bl \label{N tor}
Let $M$ be a finitely generated torsion $\Zp\ps{G}$-module, where $G=\Zp^d$. Suppose that $N$ is a closed subgroup of $G$ satisfying $N\cong\Zp$ and $G/N\cong \Zp^{d-1}$, with the additional condition that $M_N$ is torsion over $\Zp\ps{G/N}$. If we let $\ga_N$ be a topological generator of $N$, then the submodule $M[\ga_N-1]$ consisting of elements annihilated by $\ga_N-1$ is a pseudo-null $\Zp\ps{G}$-module.
\el

\bpf
By \cite[Lemma 4.5]{LimFine}, we have
\[ \rank_{\Zp\ps{G}}(M) = \rank_{\Zp\ps{G/N}}(M_N) - \rank_{\Zp\ps{G/N}}H_1(N,M).  \]
Taking the hypotheses of the lemma into account, we see that $H_1(N,M)$ is a torsion $\Zp\ps{G/N}$-module. Moreover, since $N$ is procyclic with generator $\ga_N$, we have the identification $H_1(N.M)= M^N= M[\ga_N-1]$. Now, let $H$ be a closed subgroup of $G$ such that $H$ is mapped isomorphically to $G/N$ through the natural surjection $G\twoheadrightarrow G/N$. Under this identification, the module $M[\ga_N-1]$ may be viewed as a finitely generated torsion $\Zp\ps{H}$-module. It then follows from \cite[Lemma 5.1]{LimFine} that $M[\ga_N-1]$ is a pseudo-null $\Zp\ps{G}$-module.
\epf

Building on the preceding lemma, we have the following.

\bl \label{N inject}
Retain the setting of Lemma \ref{N tor}, and let $M^0$ denote the maximal pseudo-null submodule of $M$. Then there is an injection
\[ (M^0)_N\hookrightarrow M_N.\]
\el

\bpf
Let $\ga$ be a topological generator of the group $N$. From the short exact sequence
\[ 0\lra M^0\lra M \lra M/M^0\lra 0,\]
we obtain the exact sequence
\[
  0\lra M^0[\ga-1]\lra M[\ga-1] \lra M/M^0[\ga-1]
  \lra (M^0)_{N} \lra M_{N}.\]
By Lemma \ref{N tor}, the module $M/M^0[\ga-1]$ is pseudo-null over $\Zp\ps{G}$. On the other hand, since
$M^0$ is defined as the maximal pseudo-null submodule of $M$, the quotient $M/M^0$ has no nonzero pseudo-null submodules. The combination of these two observations then forces $M/M^0[\ga-1]=0$. With this term eliminated, the exact sequence collapses to give the injection we set out to prove.
\epf

The next lemma is concerned with the descent of pseudo-null modules, when $G=\Zp^d$ with $d\geq 3$.

\bl \label{N pseudo}
Let $M$ be a finitely generated pseudo-null $\Zp\ps{G}$-module, where $G=\Zp^d$ with $d\geq 3$. Then there exist a closed subgroup $N$ of $G$ satisfying $N\cong\Zp$ and $G/N\cong \Zp^{d-1}$, with the additional condition that $M_N$ is pseudo-null over $\Zp\ps{G/N}$.
\el

\bpf
While this result is essentially well-known, we include a brief sketch of its proof here for the reader's convenience. Since $M$ is pseudo-null, it has a non-zero annihilator which does not lie in the ideal $p\Zp\ps{G}$. Therefore, we may apply \cite[Lemma 2]{Gr} to see that there exist a closed subgroup $H$ of $G$ satisfying $G/H\cong \Zp$, and with the additional condition that $M$ is finitely generated over $\Zp\ps{H}$. Subsequently, applying \cite[Lemma 2.6]{Fuj}, we obtain a closed subgroup $V$ of $G$ such that $V\subseteq H$ and $H/V\cong \Zp$, and such that $M_V$ is  pseudo-null over $\Zp\ps{G/V}$. Let $N$ be a closed subgroup of $V$ such that $N\cong\Zp$ and $G/N\cong \Zp^{d-1}$. Since $M_V= (M_N)_{V/N}$, it follows that $M_N$ is pseudo-null over $\Zp\ps{G/N}$.
\epf

For the next lemma, we specialize the case where $G=\Zp^2$. Recall that for a given $\Zp\ps{G}$-module $M$, we denote by $M^0$ the maximal pseudo-null submodule of $M$.

\bl \label{M^0_N}
Let $M$ be a finitely generated $\Zp\ps{G}$-module, where $G=\Zp^2$. Suppose that $M$ satisfies the property that $M_U$ is finite for every open subgroup of $G$. Then there exists  a closed subgroup $N$ of $G$ satisfying $N\cong\Zp$ and $G/N\cong \Zp$, with the additional condition that $(M^0)_N$ is finite.
\el

\bpf
The conclusion of the lemma is trivially true if $M^0=0$. We may therefore assume that $M^0\neq 0$. Fix a topological basis $\sigma, \tau$ of $G$ and identify $\Zp\ps{G}\cong\Zp\ps{U,W}$ via the assignment $\sigma\mapsto U$, $\tau\mapsto W$. For $\alpha\in\Zp$, let $N_\alpha$ denote the closed subgroup topologically generated by $\sigma^{-\alpha}\tau$ which in turn corresponds to the element $W_{\alpha}:=(1+U)^{-\alpha}(1+W)-1$. We shall now show that the characteristic element of $(M^0)_{N_\alpha}$ is coprime to $(1+U)^{p^m}-1$ for every $m$. The hypothesis of the lemma guarantees that $M_{N_\alpha}/(1+U)^{p^m}-1$ is finite for every $m\geq 0$. Therefore, the characteristic element of $M_{N_\alpha}$ is coprime to $(1+U)^{p^m}-1$ for every $m$. By virtue of the injection $(M^0)_{N_\alpha} \hookrightarrow M_{N_\alpha}$ as asserted by Lemma \ref{N inject}, the same coprimality assertion holds for the characteristic element of $(M^0)_{N_\alpha}$.

Now, let $0=\cap_{i=1}^r M^0/M_i$ be a shortest primary decomposition of $0$ in the Noetherian module $M^0$ (in the sense of \cite[\S 6]{Matsu}). In other words, we have each $M_i$ being submodule of $M^0$ with $M^0/M_i$ having exactly one associated prime $P_i$, and that $P_i\neq P_j$ whenever $i\neq j$. By the proof of \cite[Lemma 2]{Oz},
the natural map
\[ M^0\lra \bigoplus_{i=1}^r M^0/M_i \]
is injective and has a finite cokernel, which we shall denote by $D$. Now suppose that $N$ is a closed subgroup of $G$ satisfying $N\cong\Zp$ and $G/N\cong \Zp$. We then have an exact sequence
\[ D^N\lra (M^0)_N\lra \bigoplus_{i=1}^r \big(M^0/M_i\big)_N\lra  D_N\lra 0. \]
The question is therefore reduced to finding an appropriate subgroup $N$ such that every $(M^0/M_i\big)_N$ is finite. In fact, we will show for each $i$,
the module $(M^0/M_i\big)_{N_\alpha}$ can be infinite for at most one $\alpha$. Suppose on the contrary that $(M^0)_{N_\alpha}$ and $(M^0)_{N_\beta}$ are infinite for $\alpha\neq \beta$. Then by \cite[Lemma 3]{Oz}, one has $(1+U)^{-\alpha}(1+W)-1, (1+U)^{-\beta}(1+W)-1 \in P_i$. Furthermore, the proof of the said lemma tells us that $P_i$ is a prime ideal of height $2$ generated by $W_\alpha$ and an irreducible element $g_\alpha$ in $\Zp\ps{U}$, where $g_\alpha$ is a divisor of the characteristic element $\mathrm{char}_{\Zp\ps{S}}(M^0)_{N_\alpha}$. But we also have
\[ (1+U)^{-\beta}(1+W)\big((1+U)^{\beta-\alpha}-1\big) = (1+U)^{-\alpha}(1+W)-1 - \big((1+U)^{-\beta}(1+W)-1\big) \in P_i \]
which in turn implies that $(1+U)^{p^n}-1\in P_i$ for some integer $n\geq 0$. As $P_i= (W_\alpha, g_\alpha)$, we have $g_\alpha$ dividing $(1+U)^{p^n}-1$. However, by the argument in the opening paragraph, the element $(1+U)^{p^n}-1$ is coprime to the characteristic element of $(M^0)_{N_\alpha}$ by the argument in the initial paragraph. This forces $g_\alpha$ to be a unit, contradicting the fact that $g_\alpha$ is irreducible. This contradiction completes the proof.
 \epf

We can now give the proof of Theorem \ref{mainthm2a}.

\bpf[Proof of Theorem \ref{mainthm2a}]
Although the theorem is stated for $d\geq 2$, we find it convenient to establish the result for all $d\geq 1$, where the proof proceeds by induction. The base case $d=1$ is an immediate consequence of Theorem \ref{mainthm1}. Now suppose $d\geq 2$. By either Lemma \ref{M^0_N} or Lemma \ref{N pseudo}, depending on whether we are in the case $d=2$ or $d\geq 3$ respectively, there exists a closed subgroup $N$ of $G$ with $N\cong\Zp$ and $G/N\cong \Zp^{d-1}$, such that $(X_{2i,F_\infty}^0)_N$ is a pseudo-null $\Zp\ps{G/N}$-module.
By our underlying hypothesis, we have $X^0_{2i,F_\infty}\neq 0$, so an application of Nakayama lemma immediately implies that $(X^0_{2i,F_\infty})_N\neq 0$. Now, let $L_\infty$ denote the fixed field of $N$. It then follows from Proposition \ref{gal descent} that $(X_{2i,F_\infty})_N\cong X_{2i,L_\infty}$, and the latter is torsion over $\Zp\ps{G/N}$ by Proposition \ref{tor}. Applying Lemma \ref{N inject}, we obtain an injection
 \[ (X^0_{2i,F_\infty})_N\hookrightarrow  X_{2i,L_\infty}. \]
Consequently, we conclude that $X_{2i,L_\infty}^0\neq 0$. Applying our induction hypothesis to the $\Zp^{d-1}$-extension $L_\infty$, we obtain a finite extension $L$ of $F$ contained in $L_\infty$ with
\[ \ker\Big(K_{2i}(\Op_L)[p^\infty] \lra A_{2i, L_\infty}\Big)\neq 0.\]
 As this kernel is contained in $\ker\Big(K_{2i}(\Op_L)[p^\infty] \lra A_{2i, F_\infty}\Big)$, the conclusion of the theorem follows.
\epf

We now prove the partial converse of Theorem \ref{mainthm2a} which is also an algebraic $K$-theoretical analogue of \cite[Theorem 2]{Fuj}.

\bt \label{mainthm2b}
Let $F_\infty$ be a $\Zp^d$-extension of $F$ with $d\geq 2$. Suppose that $K_{2i}(\Op_F)[p^\infty]$ is a nontrivial cyclic group and that there exists a finite extension $L$ of $F$ contained in $F_\infty$ such that
\[ \ker\Big(K_{2i}(\Op_L)[p^\infty] \lra A_{2i, F_\infty}\Big)\neq 0.\]
Then $X^{0}_{2i,F_\infty}\neq 0$.
\et

Before proceeding to the proof of the theorem stated above, we first establish the following, where we note that this said result does not require the cyclicity assumption on $K_{2i}(\Op_F)[p^\infty]$.

\bt \label{projdim}
Let $F_\infty$ be a $\Zp^d$-extension of $F$ with $d\geq 2$. Suppose that there exists a finite extension $L$ of $F$ contained in $F_\infty$ such that
\[ \ker\Big(K_{2i}(\Op_L)[p^\infty] \lra A_{2i, F_\infty}\Big)\neq 0.\]
Then the projective dimension of $X_{2i,F_\infty}$ as a $\Zp\ps{G}$-module is at least $2$.
\et

\bpf
We begin with a preliminary observation. Since $F_\infty/F$ is a $\Zp^d$-extension, the structure of this Galois extension guarantees that we can always find a $\Zp$-extension $L_\infty$ of $L$ contained in $F_\infty$ such that
\[ \ker\Big(K_{2i}(\Op_L)[p^\infty] \lra A_{2i, L_\infty}\Big)\neq 0.\]
By Theorem \ref{mainthm1}, this implies that $X^{0}_{2i,L_\infty}\neq 0$. Taking this into account, we may apply \cite[Proposition 5.3.19]{NSW} to conclude that the projective dimension of $X_{2i,L_\infty}$ as a $\Zp\ps{\Gal(L_\infty/L)}$-module is at least $2$.

We now proceed to the proof of the theorem via contradiction. Suppose that the projective dimension of $X_{2i,F_\infty}$ as a $\Zp\ps{G}$-module is at most $1$. Since $X_{2i,F_\infty}$  is finitely generated torsion over $\Zp\ps{G}$, it admits a resolution of the form
\[ 0\lra \Zp\ps{G}^r\lra \Zp\ps{G}^r \lra X_{2i,F_\infty} \lra 0. \]
Taking $\Gal(F_\infty/L_\infty)$-coinvariant, we obtain the induced exact sequence
\[ \Zp\ps{\Gal(L_\infty/F)}^r\stackrel{f}{\lra} \Zp\ps{\Gal(L_\infty/F)}^r \lra X_{2i,L_\infty} \lra 0. \]
Since $X_{2i,L_\infty}$ is torsion as a $\Zp\ps{\Gal(L_\infty/L)}$-module, so is $\ker f$. However, $\Zp\ps{\Gal(L_\infty/F)}\cong \Zp\ps{\Gal(L_\infty/L)}^{|L:F|}$ is $\Zp\ps{\Gal(L_\infty/L)}$-torsionfree. As a result, we have $\ker f=0$, and the above sequence becomes short exact. Thus, the projective dimension of $X_{2i,L_\infty}$ as a $\Zp\ps{\Gal(L_\infty/L)}$-module is at most $1$, contradicting the earlier observation in the initial paragraph. This contradiction completes the proof.
\epf

We briefly deviate from the current discussion to make a few remarks.

\br \label{H Iw remark}
\begin{enumerate}
  \item[$(1)$] If $\Ga\cong\Zp$, and $M$ is a finitely generated $\Zp\ps{\Ga}$-module, then it is well-known that $M^0=0$ if and only if the $\mathrm{pd}_{\Zp\ps{\Ga}}(M)\leq 1$ (for instances, see \cite[Proposition 5.3.19]{NSW}). This characterization no longer holds when we move to higher dimension: for $G=\Zp^d$ with $d\geq 2$, we can only guarantee that projective dimension $\leq 1$ forces $M^0=0$, but the converse direction may fail to hold (see \cite[Lemma 4.1]{Kl}). For this reason, the non-vanishing assertion $X_{2i,F_\infty}^0\neq $ does not follow as an immediate consequence of the claims of Theorem \ref{projdim}, and thus cannot be directly inferred from it.
      \item[$(2)$] Although this remark falls outside the paper's central theme, we still include it here, as we thought it of interest. As seen in the discussion in Section \ref{K and Galois}, the module $X_{2i,F_\infty}$ can be identified as a second Iwasawa cohomology group $H^2_\Iw(F_\infty/F,\Zp(i+1))$. In Iwasawa theory, the second Iwasawa cohomology group of a Galois representation $T$ is an important object of study, for instances, it appears in one formulation of the main conjecture (see \cite[Conjecture 12.10]{Ka}). On the other hand, this group is closely related to the torsion submodule of the Pontryagin dual of a certain discrete Galois cohomology group (see \cite[Proposition 3.5]{OV}) which appears in the definition of a Selmer group. Despite its importance, the general structure of these Iwasawa cohomology groups remains poorly understood, beyond the well-known fact that their torsion property is equivalent to the weak Leopoldt conjecture for $T$ (see \cite[Theorem 3.3]{OV} or \cite[Lemma 7.1]{LimFine}). In this special case $T=\Zp(i+1)$, the conclusion of Theorem \ref{projdim} appears to be the first instance in the literature of a structural observation for the second Iwasawa cohomology group that goes beyond previously known results.
\end{enumerate}
\er

We return to our initial discussion, and finish the proof of Theorem \ref{mainthm2b}.

\bpf[Proof of Theorem \ref{mainthm2b}]
We prove this by contradiction. In particular, we assume that $X_{2i,F_\infty}$ has no nonzero pseudo-null $\Zp\ps{G}$-submodules.
Since $(X_{2i,F_\infty})_G\cong K_{2i}(\Op_F)[p^\infty]$ is cyclic, it follows from Nakayama lemma that $X_{2i,F_\infty}$ itself is a cyclic $\Zp\ps{G}$-module. Therefore, we have a short exact sequence
\[ 0\lra I\lra\Zp\ps{G} \lra X_{2i,F_\infty}\lra 0, \]
where $I$ is an ideal of $\Zp\ps{G}$. Fix a set of generators $f_1,...,f_k$ of the ideal $I$. Let $f$ be a greatest common divisor of $f_1,...,f_k$, and let $I_0$ be the ideal generated by $f_1/f,...,f_k/f$. A straightforward verification shows that the map
\[ \Zp\ps{G} \lra f\Zp\ps{G}/I,  x\mapsto xf + I\]
is surjective with kernel $I_0$. Since the elements $f_1/f,...,f_k/f$ have no common divisor, the ideal $I_0$ is not contained in any height $1$ prime. Therefore, $\Zp\ps{G}/I_0\cong f\Zp\ps{G}/I$ is a pseudo-null $\Zp\ps{G}$-submodule of $\Zp\ps{G}/I=X_{2i,F_\infty}$. Our initial assumption forces this pseudo-null submodule to be zero, which in turn implies that $I$ is a principal ideal generated by $f$. In particular, $I$ is a free $\Zp\ps{G}$-module and and therefore the projective dimension of $X_{2i,F_\infty}$ is $\leq 1$. This is a contradiction to Theorem \ref{projdim}.
\epf

\section{Application to class groups} \label{sec:class}
In this section, we present an application of our Theorem \ref{mainthm2b} to derive a sufficient criterion for the existence of non-trivial pseudo-null submodules in the classical Iwasawa module.

For a number field $F$, we let $A'_F$ denote the $p$-primary part of the ideal class group of $\Op_F[1/p]$. Write $E=F(\mu_p)$, and set $r=|E:F|$. Throughout our subsequent discussion, we maintain the standing assumption that $r>1$. Note that since $p$ is odd, we always have $E\cap F_\infty = F$ for a given $\Zp^d$-extension $F_\infty$ of $F$.
We then define $E_\infty = F_\infty(\mu_p)$. For every finite extension $L$ of $F$ contained in $F_\infty$, we identify $\Gal(L(\mu_p)/L)$ with $\Gal(F(\mu_p)/F)$. By abuse of notation, we shall denote these groups as $\Delta$.

The modules $A'_{E_\infty}$ and $X'_{E_\infty}$ are then similarly defined which come equipped with $\Zp\ps{\Gal(E_\infty/F)}$-module structures. Depending on the specific context of our discussion, we sometimes also regard these as $\Zp\ps{\Gal(E_\infty/E)}$-modules.

Once and for all, we fix a primitive $p$-th root of unity $\zeta$ in $\overline{\Q}$, and fix a generator $g$ of $(\Z/p\Z)^{\times}$. We then write $\sigma$ for the corresponding element in $\Gal(\Q(\mu_p)/\Q)$ which sends $\zeta$ to $\zeta^g$. Under this setup, the elements of the group $\Delta$ may be identified with elements of the form $\sigma^{k(p-1)/r}$, where $k=0,1,..., r-1$. Correspondingly, we have $g^{(p-1)/r}$ being a generator of the subgroup of $(\Z/p\Z)^{\times}$, which is the image of $\Delta$ under the embedding of $\Delta$ into $(\Z/p\Z)^{\times}$.

Let $\omega$ be the Teichm\"{u}ller character of the group $(\Z/p\Z)^{\times}$. For a given integer $j$, we define
\[ \e_j:=\e_j(F): = \frac{1}{r}\sum_{k=0}^{r-1}\omega(g)^{jk(p-1)/r}\sigma^{-k(p-1)/r}\]
which lives in $\Zp[\Gal(F(\mu_{p})/F)]$. One immediately sees that $\e_j=\e_{j'}$ whenever $j=j'$ (mod $r$). Furthermore, it can be checked easily that $\e_0, \e_1, ...,\e_{r-1}$ forms a collection of primitive idempotents of the group ring $\Zp[\Gal(F(\mu_{p})/F)]$.

In what follows, we will collect certain results on the $p^n$-rank of $K_{2i}(\Op_F)$ for a given number field $F$ and integer $i\geq 1$. We define $a:=a(F)$ to be the largest integer such that $F(\mu_p) = F(\mu_{p^{a}})$. For a nonzero integer $j$, denote by $S_p^{(j)}(F)$ the set of primes $v$ of $F$ above $p$ such that $j$ is divisible by $|\Delta_v|$, with $\Delta_v$ being the decomposition group of $\Gal(F(\mu_p)/F)$ at $v$. With these definitions in hand, we can now state the following.

\bp \label{p^n rank thm for small n} Suppose that $n\leq a(F)$.
\begin{itemize}
  \item[$(i)$] If $i\equiv 0$ mod $|F(\mu_p):F|$, we have the following exact sequence
\[ 0\lra \mu_{p^n}^{\ot i}\ot A'_F\lra  K_{2i}(\Op_F)/p^n \lra \bigoplus_{v\in S_p}\mu_{p^n}^{\ot i} \lra \mu_{p^n}^{\ot i} \lra 0.\]
  \item[$(ii)$] If $i\not\equiv 0$ mod $|F(\mu_p):F|$, we have the following exact sequence
\[ 0\lra \mu_{p^n}^{\ot i}\ot \e_{-i}A'_{F(\mu_p)} \lra  K_{2i}(\Op_F)/p^n \lra \bigoplus_{v\in S_p^{(i)}}\mu_{p^n}^{\ot i} \lra 0.\]
\end{itemize}
\ep

\bpf
For $i = 1$, this is a classical result of Keune \cite[Theorem 5.4]{Ke}. For the higher $K$-groups, this is established
by the author in \cite[Theorem 4.2]{LimKclass}.
\epf

\begin{itemize}
  \item[$(\mathbf{Sp})$] Every prime of $F(\ze+\ze^{-1})$ above $p$ does not split in $F(\mu_p)$.
\end{itemize}

Note that the validity of $(\mathbf{Sp})$ also implies that $r=|F(\mu_p):F|$ is even. The proof of the next lemma is straightforward (recall that the prime $p$ is always odd by our standing assumption) and is left to the reader.

\bl \label{Sp L/F}
Let $L$ be a finite Galois $p$-extension of $F$. If $F$ satisfies $(\mathbf{Sp})$, so does $L$.
\el

For our purposes, we also need the following which identifies the odd eigenspace of $A_{F(\mu_p)}$ and $A'_{F(\mu_p)}$.

\bp \label{kummer lemma}
Assume that $(\mathbf{Sp})$ is valid. Then for every odd integer $i$ such that $1\leq i\leq r-1$, we have
\[ \e_{-i}A'_{F(\mu_p)}= \e_{-i}A_{F(\mu_p)}\]
\ep

\bpf
Since $r=|F(\mu_p):F|$ is even, the proposition is equivalent to saying that $\e_{i}A'_E= \e_{i}A_E$ for every odd $i$.
The natural surjection $\pi:\Cl(\Op_{F(\mu_p)})\twoheadrightarrow \Cl(\Op_{F(\mu_p)}\big[\frac{1}{p}\big])$ has kernel which is generated by the classes of primes of $F(\mu_p)$ above $p$. Now, for an odd integer $i$, we have the following computation

\begin{align*} \label{e_i}
 \e_{i} = & ~\frac{1}{r}\sum_{k=0}^{r-1}\om(g)^{-ik(p-1)/r}\sigma^{-k(p-1)/r} \\
   = & ~\frac{1}{r}\sum_{k=0}^{\frac{r}{2}-1}\left(\om(g)^{ik(p-1)/r}\sigma^{-k(p-1)/d} + \om(g)^{i(k+\frac{r}{2})(p-1)/r}\sigma^{-(k+\frac{r}{2})(p-1)/r}\right) \\
 =& ~\frac{1}{r}\sum_{k=0}^{\frac{r}{2}-1}\om(g)^{ik(p-1)/r}\sigma^{-k(p-1)/r}\left( 1 + (-1)^i\sigma^{-\frac{p-1}{2}}\right) \\
=& ~\frac{1}{r}\sum_{k=0}^{\frac{r}{2}-1}\om(g)^{ik(p-1)/d}\sigma^{-k(p-1)/r}\left( 1 -\sigma^{\frac{p-1}{2}}\right),
\end{align*}
where we note that $\sigma^{-\frac{p-1}{2}} = \sigma^{\frac{p-1}{2}}$ is the generator of the group $\Gal(F(\ze_p)/F(\ze_p+\ze_p^{-1}))$.
By our hypothesis $(\mathbf{Sp})$, the primes of $F(\ze +\ze^{-1})$ above $p$ do not split in $F(\ze_p)$. As such, they are invariant under the automorphism $\sigma^{\frac{p-1}{2}}$, and it follows from this that $\ker \pi$ is annihilated by $1 -\sigma^{\frac{p-1}{2}}$. As a consequence, $(\ker \pi)[p^\infty]$ is annihilated by $\e_i$ for every odd $i$, which in turn yields the required identification $\e_{i}A'_E= \e_{i}A_E$.
\epf

We can now prove the following variant of \cite[Theorem 2]{Fuj}.

\bt \label{class group}
Assume that $(\mathbf{Sp})$ is valid for $F$. Let $F_\infty$ be a $\Zp^d$-extension of $F$ such that $F^\cyc\subseteq F_\infty$. Suppose that there exists a positive odd integer $i$ with $1\leq i\leq r-1$ that satisfies all of the following conditions.
 \begin{enumerate}
   \item[$(a)$] There exists a finite extension $L$ of $F$ contained in $F_\infty$ such that
      \[ \ker\Big(K_{2i}(\Op_{L})[p^\infty] \lra A_{2i,F_\infty}\Big)\neq 0.\]
\item[$(b)$] $K_{2i}(\Op_F)[p^\infty]$ is a non-trivial cyclic group.
  \end{enumerate}
Then $X_{E_\infty}^0\neq 0$.
\et

\bpf
From the natural decomposition
\[ X_{E_\infty} = \bigoplus_{k=0}^{r-1}\e_k X_{E_\infty}, \]
 we are therefore reduced to showing that $\e_{-i} X_{E_\infty}=\e_{r-i} X_{E_\infty}$ has a non-trivial pseudo-null submodule.

For every finite extension $L$ of $F$ contained in $F_\infty$, Proposition \ref{p^n rank thm for small n} gives us the following short exact sequence
 \[ 0\lra \mu_{p^{n(L)}}^{\ot i}\ot \e_{-i}A_{L(\mu_p)} \lra  K_{2i}(\Op_L)/p^{n(L)} \lra \bigoplus_{v\in S_p^{(i)}(L)}\mu_{p^{n(L)}}^{\ot i} \lra 0,\]
 where we have identified $\e_{-i}A'_{L(\mu_p)}$ with $\e_{-i}A_{L(\mu_p)}$ by combining Lemma \ref{Sp L/F} and Proposition \ref{kummer lemma}. Also, as a consequence of $L$ satisfying $(\mathbf{Sp})$, the decomposition group of $\Gal(L(\mu_p)/L)$ at a prime above $p$ has even order. Since $i$ is odd, this in turn implies that the set $S_p^{(i)}(L)$ is empty, and so we have an identification
 \begin{equation}\label{eigenspace K}
   \mu_{p^{a(L)}}^{\ot i}\ot \e_{-i}A_{L(\mu_p)} \cong  K_{2i}(\Op_L)/p^{a(L)}.
 \end{equation}
Since $F_\infty$ contains $F^\cyc$, the field $E_\infty$ contains all $p$-power root of unity, and so $a(L)\to \infty$. Thus, upon taking inverse limit, we obtain
 \[ \big(\e_{-i}X_{E_\infty}\big)(i) \cong X_{2i, F_\infty}. \]
 Clearly, $(\e_{-i}X_{E_\infty})(i)^0\neq 0$ if and only if $(\e_{-i}X_{E_\infty})^0\neq 0$. Combining this latter observation with the preceding isomorphism, we are thus reduced to proving that $X_{2i, F_\infty}^0\neq 0$. Under the hypotheses of the present theorem, this nonvanishing property follows directly as a consequence of Theorem \ref{mainthm2b}.
\epf

We end with the following remarks.

\br
\begin{itemize}
     \item[$(1)$] It is evident from the proof of Theorem \ref{class group} that we have actually showed that  $(\e_{-i}X_{E_\infty})^0\neq 0$.
     \item[$(2)$] In view of (\ref{eigenspace K}), the condition that $K_{2i}(\Op_F)[p^\infty]$ is a non-trivial cyclic group is equivalent to the statement that $\e_{-i}A_{E}$ itself is a non-trivial cyclic group. This is hence a weaker assumption than having the full group $A_E$ being non-trivial cyclic.
     \item[$(3)$] By the preceding remark, it would be of interest to ask whether the hypothesis of the theorem can be phrased entirely on capitulation of ideal class groups. This brings us to the question whether the capitulation on $K$-groups is related with the capitulation of ideal class groups. However, as seen in (\ref{eigenspace K}), at a finite level, what we can obtain is only an identification of a twist of an appropriate eigenspace of $A_{L(\mu_p)}$ with $K_{2i}(\Op_L)[p^\infty]$ modulo some $p$-power. This partial identification, while useful, fails to capture the full group structure, and for that reason we are currently unable to establish the exact relationship between the two capitulation phenomena. As a result, the initial problem we posed stays open, and we can only hope to come back to this in a subsequent work.
   \end{itemize}
\er

\end{document}